\documentclass[conference]{IEEEtran}
\IEEEoverridecommandlockouts
\PassOptionsToPackage{hyphens}{url}
\usepackage[breaklinks]{hyperref}
\usepackage{tikz}
\usepackage{pgfplots}
\usepackage{cite}
\usepackage{amsmath,amssymb,amsfonts}
\usepackage{algorithm}
\usepackage{algorithmic}
\usepackage{graphicx}
\usepackage{textcomp}
\usepackage{xcolor}
\usepackage{physics}
\def\BibTeX{{\rm B\kern-.05em{\sc i\kern-.025em b}\kern-.08em
    T\kern-.1667em\lower.7ex\hbox{E}\kern-.125emX}}

\newcommand{\binary}{\{0, 1\}}

\newcommand{\secref}[1]{Section~\ref{#1}}
\newcommand{\fig}[1]{Figure~\ref{#1}}
\newcommand{\tab}[1]{Table~\ref{#1}}
\newcommand{\cm}{\textrm{cm}}
\newcommand{\constraintsButX}{\Gamma(y,z)}
\newcommand{\fixedBoxes}{\mathcal F}

\begin{document}

\title{Accelerated Benders Decomposition for\\Variable-Height Transport Packaging Optimisation}

\author{
    \IEEEauthorblockN{Alain Lehmann\textsuperscript{\textsection}, Wilhelm Kleiminger}
    \IEEEauthorblockA{\textit{Ergon Informatik AG} \\
    Merkurstrasse 43, CH-8032 Zurich\\
    \{alain.lehmann, wilhelm.kleiminger\}@ergon.ch
    }
    \and
    \IEEEauthorblockN{Hakim Invernizzi, Aurel Gautschi}
    \IEEEauthorblockA{\textit{Digitec Galaxus AG} \\
    Pfingstweidstrasse 60b, CH-8005 Zürich\\
    \{Hakim.Invernizzi, Aurel.Gautschi\}@digitecgalaxus.ch
    }
}

\maketitle

\begingroup\renewcommand\thefootnote{\textsection}
\footnotetext{Main author}
\endgroup

\begin{abstract}
    This paper tackles the problem of finding optimal variable-height transport
packaging. The goal is to reduce the empty space left in a box when shipping
goods to customers, thereby saving on filler and reducing waste.
We cast this problem as a large-scale mixed integer problem (with over seven
billion variables) and demonstrate various acceleration techniques to solve it
efficiently in about three hours on a laptop.
We present a KD-Tree algorithm to avoid exhaustive grid evaluation of the 
3D-bin-packing, provide analytical transformations to accelerate the Benders
decomposition, and an efficient implementation of the Benders sub problem for
significant memory savings and a three order of magnitude runtime speedup.

\end{abstract}


\section{Introduction}

Online retailers 
use transport packaging to safely ship products from their logistics centres to their customers. The right choice of transport packaging ranges from none at all (if the product is already well protected by its factory packaging) over padded envelopes to boxes made from heavy carton stuffed with filler material for additional protection~\cite{Gurumoorthy2020}. 

Choosing the right packaging is important as the environmental footprint of the transport packaging alone can range from 4.1\% to 26\% of the order process~\cite{buycom,VANLOON2015478}. There is also an economic incentive -- lower shipping costs~\cite{ortec}. 
The ratio of product to air in a box is a major factor to account for.
More air means more packaging and filler material. Thus a greater environmental footprint and higher operational costs. 

However, having a myriad of box sizes to fit each combination of products is impractical as the space at the packers' workstations is limited.
Instead, a fixed set of variable-height boxes is used with carton dimensions optimised for the expected deliveries\footnote{Automated packing system that fit transport packing to products are not applicable for bulky goods and scalable during peak season.}.

Optimising variable-height transport packaging can be formulated as a mixed integer program (MIP) similar to the facility location problem.
The caveat is that the number of possible box sizes and the expected deliveries sample make this
a very large-scale optimisation problem. We show various acceleration techniques to make 
the problem tractable. Specifically, we make the following contributions (c.f. \tab{tab:contributions}):

\begin{enumerate}
    \item model variable-height transport packaging as a mixed integer problem,
    \item accelerate its Benders decomposition using analytic transformations for a 42-fold speedup,
    \item implement the Benders sub program efficiently with 1000-fold speedup and 64-fold memory reduction,
    \item propose an adaptive KD-Tree algorithm to bin pack every combination of packing unit to box with a 29 fold speedup.
\end{enumerate}

\begin{figure}[tb]
    \centering
    \includegraphics[width=\columnwidth]{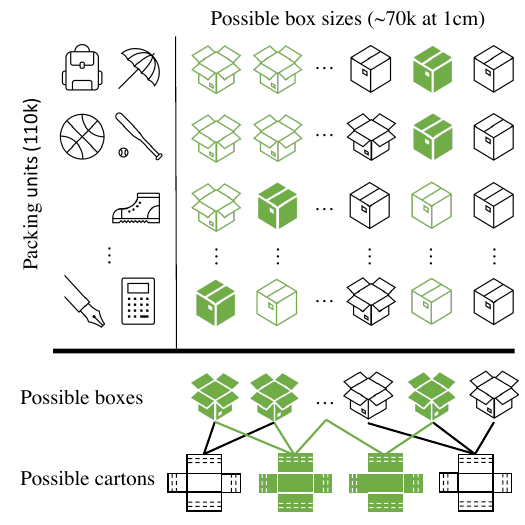}%
    \caption{
        Overview of the carton box optimisation model (c.f. \secref{sec:problemstatement}).
        The goal is to find optimal variable-height boxes to transport packing units
        with minimal empty space.
        (bottom): Illustration of cartons and boxes with ``selected'' ones coloured green.
        Lines connecting cartons and boxes express the variable-height box relationship.
        (top-right): Illustration of the ``fitting matrix'' $F$ (shape)
        and the packing variables $x$ (color): Green coloured boxes indicate that a box $b$
        is available for packing ($y_b=1$), while black ones are not available. Closed
        boxes indicate packing unit $p$ does fit into a box ($F_{pb}=1$), while open ones do not.
        Closed, dark green boxes are assumed for shipping $x_{pb}=1$ and contribute to
        the overall objective value.
    }
    \label{fig:overview}
\end{figure}

\begin{table}[htbp]
    \caption{Contribution Summary of Performance Improvements}
\begin{center}
\begin{tabular}{|l|c|r|r|r|}
\hline
    Improvement & Section &Before & After & Speedup \\
\hline
    Bender Overall Optimiation & \ref{sec:bender} & $8356s$ & $199s$ & 42 \\
    - Runtime Sub Problem & \ref{sec:subprogram:runtime} & 7s & 0.007s & 1000 \\
    - Memory Sub Problem & \ref{sec:subprogram:memory} & 60GB & $<1$GB & 64 \\
\hline
    KD-Tree Bin Packing      & \ref{sec:kdtree} & & &  \\
    - CPU-Time (i.e. 1-core) & & 535.76h & 18.40h & 29 \\
    - On M1 Pro/10-cores     & & 74h17m & 2h59 & 25 \\
\hline
\end{tabular}
\end{center}
\label{tab:contributions}
\end{table}


This paper is structured as follows: \secref{sec:problemstatement} presents
our model and discusses related work. \secref{sec:bender} revisits
the Benders decomposition for our concrete model and shows how to handle
variable-height boxes efficiently. \secref{sec:subprogram} derives the analytical solution
of the Benders sub program and outlines an efficient implementation.
\secref{sec:kdtree} discusses an adaptive KD-Tree algorithm to compute the 
3D-bin-packing matrix. Finally, \secref{sec:results} evaluates the proposed 
contributions on a large scale dataset, followed by future work and conclusion
in \secref{sec:conclusion}.

\section{Problem Statement: Carton Box Optimisation}\label{sec:problemstatement}

\begin{figure}[tb]
    \centering
    \includegraphics[width=0.35\columnwidth]{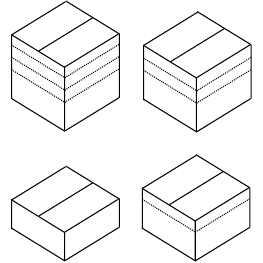}%
    \hfill
    \begin{tikzpicture}
        \small
        \begin{axis}[
                name=chart,
                height=0.4\columnwidth,
                width=0.5\columnwidth,
                ytick={5000,10000,...,25000},
                ybar,
                ymin=0,
                xtick=data,
                scaled y ticks=false,
                bar width=12pt,
                enlarge x limits=0.2,
                xlabel={Number of boxes per carton},
                ylabel={\# cartons}
            ]
            \addplot[blue, fill=blue!10] plot coordinates {
                (1, 921)
                (2, 11031)
                (3, 25412)
                (4, 5892)
            };
        \end{axis}
    \end{tikzpicture}
    \caption{%
				(left) %
        Crease lines quadruple the number of box formats possible given a single carton.
        With a few small cuts, more box formats are available without more space requirements.
				(right) %
        Statistics of cartons and boxes. The 1cm discretisation yields 71790 boxes.
        The 43256 cartons are derived from boxes: Every box yields a carton unless
        that carton (with its crease lines) is a already covered by a bigger carton.
        Encoding the carton box relationship is a table of 122787 entries.
    }
\label{fig:creaseline:histogram}
\label{fig:creaselines}
\end{figure}

This section presents the model (c.f. \fig{fig:overview}) for finding optimal
variable-height transport packing, i.e. cartons with crease lines 
(c.f. \fig{fig:creaselines}) that allow for multiple boxes given a single carton.
Overall, the goal is to find a set of cartons that minimise the empty space
when shipping a packing unit in a box built from one of the optimal cartons.

The main input is a fitting matrix $F_{pb}$ that indicates if (the items of)
packing unit $p=1\dots P$ fit into box $b=1\dots B$. This assumes a
discretisation of the set of possible boxes. Note that the actual
shape of the box and packing unit items is only relevant for computing the
boolean matrix $F_{pb}$ which will be discussed in \secref{sec:kdtree}.
For the optimisation, only the free space $C_{pb}=V_b-V_p$ left after packing
all items with combined volume $V_p$ into the box with volume $V_b$ is
relevant.

Minimising the overall empty space to reduce the need for filler, can thus be 
formulated as a mixed integer program (MIP)
\begin{subequations}\label{eq:box}
    \begin{align}
         \min_{x,y,z}      \sum_{p}\sum_b & C_{pb} x_{pb}               &\textrm{subject to}\label{eq:box:objective}\\
        \forall p     \;\quad  \qquad \sum_{b} F_{pb} x_{pb} & = 1 &\textrm{packing unit shippable}\label{eq:box:shippable} \\
        \forall(p,b) \quad\qquad\qquad x_{pb} & \leq y_b                 &\textrm{box available}\label{eq:box:available} \\
        \forall (z_k, y_b) \in REL  \quad z_k &\leq y_b         &\textrm{carton implies boxes}\label{eq:carton:c2b}\\
        \forall b  \qquad        \sum_{ (z_k, y_b) \in REL } z_k  &\geq  y_b &\textrm{box requires carton}\label{eq:carton:b2c}\\
                          \sum_k z_k & = M &\textrm{limited cartons}\label{eq:carton:limit}\\
        \forall b \in \fixedBoxes  \qquad\qquad y_b &= 1 &\textrm{fixed boxes}\label{eq:box:fixed}
    \end{align}
\end{subequations}
whose variables $x$, $y$, $z$, additional inputs $M$, $REL$, $\fixedBoxes$ and related
constraints are detailed in the sequel while a summary is provided in \tab{tab:notation}.

The \emph{box variables} $y_b\in\binary$ and \emph{carton variables}
$z_i\in\binary$ indicate whether a box $b$ and carton $i$ are selected for the
optimal solution, while the \emph{packing variable} $x_{pb}\in\binary$ defines
the box assumed for shipping packing unit $p$.

The constraints \eqref{eq:box:shippable} ensure that every packing unit is
shippable, i.e., there is exactly one non-zero packing $x_{pb}=1$ per packing
unit $p$ and its items fit into that specific box $F_{pb}=1$. This further
implies that the objective \eqref{eq:box:objective} has $P$ non-zero terms that
provide the total empty space for a given solution. Constraints
\eqref{eq:box:available} enforce that only boxes selected in the optimal
solution are used for shipping.

Constraint \eqref{eq:carton:limit} allows exactly $M$ cartons in the final
solution. Obviously, better solutions (with less empty space) can be found when
allowing more cartons. But the packers' workspace imposes physical restrictions on
what is possible.

Variable-height boxes are an important tool to better leverage the
available space. They quadruple the number of box sizes that can be stored in
the same space. We model variable-height boxes with crease lines by defining
a relationship between carton and boxes, i.e. a
table $REL=[(y_b, z_k)]$ that defines which boxes $y_b$ can be constructed from
carton $z_k$. This list of pairs is provided as input to the optimisation problem
and is derived from a custom logic. \fig{fig:creaseline:histogram} provides some
statistics about our concrete choice. Constraints \eqref{eq:carton:c2b} then
define that a box is part of the optimal solution whenever a related carton is,
while \eqref{eq:carton:b2c} ensures that at least one related carton is part of
the solution whenever a box is.

Finally, constraint \eqref{eq:box:fixed} defines explicit business requirements
to include a set $\fixedBoxes$ of certain boxes in the final solution. For example, some
boxes still fit into a ``Milchkästli''\footnote{\url{https://de.wikipedia.org/wiki/Milchkasten}} or benefit from special shipping fees. Moreover,
we require the largest box to always be selected.

Modelling the optimisation problem is straightforward, but solving
it is challenging due to the huge number of optimisation variables.
Specifically, the packing variable $x$ scales with the number of packing units
$P$ and the number of boxes $B$, which alone is about $B\times P =70k \times 110k =
7700 M$. Solving this large scale problem is the key contribution of this paper.

\begin{table}[tb]
    \caption{Notation}
\begin{center}
\begin{tabular}{|c|l|}
\hline
$B$ & Total number of boxes\\
$P$ & Total number of packing units \\
$K$ & Total number of cartons\\
\hline
$b$ & Index variable used for boxes\\
$p$ & Index variable used for packing units\\
$k, l$ & Index variables used for cartons\\
\hline
$Volume(b)$ & Volume of a box \\
$Volume(p)$ & Volume of a packing unit (sum of item volumes)\\
$F_{pb}$ & Does packing unit $p$ fit into box $b$\\
$C_{pb}$ & Cost when fitting packing unit $p$ into box $b$\\
$M$ & Maximal allowed number of cartons\\
\hline
$x_{pb}$ & Is packing unit $p$ packed into box $b$\\
$y_{b}$ & Is box $b$ selected for production\\
$z_{k}$ & Is carton $k$ selected for production\\
$REL$ & Table encoding box-carton relationship\\
\hline
$\pi_p, \mu_{pb}$ & Dual variables related to packing constraints \\
$\pi^i_p, \mu^i_{pb}$ & Dual variables for box selection $y^i$\\
$s^i, w^i$ & Bender cuts generated in iteration $i$\\
\hline
\end{tabular}
\label{tab:notation}
\end{center}
\end{table}

\subsection{Related Work}\label{sec:relatedwork}

Benders decomposition \cite{benders} is a technique for large scale mixed
integer programs that has gained increased interest in recent years
\cite{bendersSurvey}. Extensions like Multi-Cut aggregations
\cite{adaptivemulticut} have been proposed to improve efficiency in general,
while this paper shows that problem-specific optimisation can have a tremendous
performance impact.

The proposed model is closely related to the well-studied factory location
problem \cite{Fischetti17}. In that context, packing units relate to customers
and boxes are the factory locations. Our model differs due to the additional
``cartons'' that represent variable-height boxes and the ``separable'' structure of the cost
coefficients. This paper demonstrates optimisation for this specific case.

Modern MIP solvers even include automatic Benders decomposition
\cite{bendersCplex}. This is generally helpful, but it may not exploit all
the structures of a concrete problem. Specifically, analytic solutions, or
bit-packed representation for memory efficiency that we use are not to be expected.

\section{Benders decomposition}\label{sec:bender}

Benders decomposition \cite{benders} is a technique to solve large scale
mixed integer programs.
On a conception level, Benders decomposition replaces continuous variables of a
mixed integer program by a single continuous variable, at the expense of
additional constraints. This reduces the number of optimisation variables which
is key to make a problem tractable. While the set of additional constraints is
huge, most are inactive for the optimal solution. Hence, Bender proposed to
iteratively find violated constraints instead of enumerating them all.

Consequently, Benders decomposition is an iterative algorithm that
alternately solves a mixed integer master program and a continuous linear sub
program. The master program is a relaxed version of the original problem. It
approximates the objective function and generates solution candidates. The sub
program scores the candidates and provides violated
constraints to improve the master program's approximation.

This section is specialised to our concrete mixed integer program \eqref{eq:box} from
\secref{sec:problemstatement}, while we refer to the survey \cite{bendersSurvey} for 
a more general description. 
\secref{sec:bender:packing} derives the Benders decomposition that eliminates 
the packing variables $x$ from the mixed integer program. \secref{sec:bender:carton}
builds on that result and also moves the box variables $y$ into the sub program,
leaving only the carton variables $z$ in the master program.


\subsection{Benders decomposition on packing variables}\label{sec:bender:packing}

Using Benders decomposition requires a set of continuous variables, 
but the problem statement in \secref{sec:problemstatement}
assumed all variables to be binary. It's thus worth noting that the \emph{optimal}
solution of \eqref{eq:box} will have binary values for the packing variable $x$,
even if the variable is relaxed to be continuous.%
\footnote{A fractional result would imply to not put everything in
the smallest available box, which would only worsens the objective value.} %
Benders decomposition transforms the original program 
into 
\begin{subequations}\label{eq:master:nonlinear}
\begin{align}
    \min_{y, z} \quad  & f_{primal}(y) & \text{subject to} \quad \constraintsButX
\end{align}
\end{subequations}
where $\constraintsButX$ collects all constraints
(\ref{eq:carton:c2b},~\ref{eq:carton:b2c},~\ref{eq:carton:limit},~\ref{eq:box:fixed})
not involving $x$, and the objective is a linear program
\begin{subequations}
\begin{align}\label{eq:sub:primal}
    f_{primal}(y)  = \min_x \sum_{pb} c_{pb} x_{pb} 
        & \quad \text{subject to} \\
    x_{pb} \leq y_b \quad \forall (p, b) 
    \quad \textrm{and} \quad &
    \sum_{b} F_{pb} x_{pb} = 1 \quad\forall p
\end{align}
\end{subequations}
that includes all constraints (\ref{eq:box:shippable},~\ref{eq:box:available})
involving $x$.
This sub program is non-linear in $y$ making \eqref{eq:master:nonlinear} a more 
complicated non-linear integer program. However, its dual
\begin{equation}\label{eq:sub:dual}
\begin{aligned}
    f_{dual}(y) = \max_{\pi, \mu\leq 0} \sum_p \pi_p + \sum_p \sum_{b} y_b \mu_{pb}\\
\text{subject to}\quad
\mu_{pb} \leq c_{pb} - F_{pb}\pi_p \quad\forall (p,b)
\end{aligned}
\end{equation}
with dual variables $\mu\in \mathcal R^{P\times B}$ and $\pi\in \mathcal R^P$ reveals that it is convex in $y$, because
the feasible region does not depend on $y$ and the objective is a maximum of
linear functions. This allows to lower bound $f$ by hyperplanes
\begin{equation}\label{eq:sub:lowerbound}
\begin{aligned}
    f(y) \geq f(y^i) + \left(\frac{df}{dy}(y^i)\right)'(y-y^i) = s^t+ (w^i)'y \; \\
    \textrm{with} \; 
    w^i_b := \frac{df}{dy_b}(y^i) = \sum_p \mu_{pb},
    \quad s^i := f(y^i) - (w^i)'y^i\\
\end{aligned}
\end{equation}
for every feasible solution $y^i$ and $'$ denotes transposition.
Benders decomposition introduces a continuous variable $\theta\in\mathcal R$
to replace the non-linear objective to obtain the integer linear master program
\begin{equation}
\begin{aligned}\label{eq:master:linear}
    \min_{y,z,\theta} \quad & \theta \quad \text{subject to}\\
        \constraintsButX 
        \quad \textrm{and} \quad & 
        \theta \geq s + w'y \quad \forall (s,w) \in OC
\end{aligned}
\end{equation}
with optimality cuts $OC$ defined by \eqref{eq:sub:lowerbound} that bound the
objective $\theta$ from below. In general (c.f.~\cite{bendersSurvey}), the
Benders decomposition includes feasibility cuts (of the form $0\leq s + w'y$)
derived from extreme rays of infeasible dual sub problems. In our concrete
setting, the business requirement (c.f. \eqref{eq:box:fixed}) to always
include the largest box $y_{b_{\max}}=1$ implies that all sub problems
are always feasible, and no further feasibility cuts are necessary.

Noting that $s^i=f(y^i)-(w^i)'y^i = \sum_p \pi^i_p$, the
Bender master program without packing variable $x$ becomes
\begin{equation}
\begin{aligned}
    \min_{y,z,\theta}\; & \theta \quad \text{subject to}\quad\\
    & \constraintsButX  \textrm{ and } y_{b_{max}}=1 \textrm { and}\\
    & \theta \geq \sum_p \pi^i_p + \sum_{pb} y_b \mu^i_{pb} \quad\forall i
\end{aligned}
\end{equation}
with $\pi^i, \mu^i$ the optimal dual variables obtained by validating
the best master solution $y^i$ in iteration $i$ with the sub program.
This mixed integer program is much smaller than the original one \eqref{eq:box}
as the number of variables no longer scales with the number of packing units (c.f. \tab{tab:bender:summary}).
However, it still has many constraints which we will address in the next section.

\subsection{Benders decomposition on packing and box variables}\label{sec:bender:carton}

This section extends the Benders decomposition to move both the packing
variables $x$ and the box variables $y$ to the sub program. This leaves
only the carton variables $z$ in the master program and also reduces the
number of constraints. This provides significant advantages to the master
program:

\begin{itemize}
    \item Less than half the number of variables (there are more boxes than cartons)
    \item No more carton-box relationship constraints (\ref{eq:carton:c2b},~\ref{eq:carton:b2c})
        which eliminates $\order{B+K}$ constraints.
\end{itemize}

Remember that moving variables to the sub problem requires them to be
continuous. The ``bi-directional'' implication defined by the crease line
constraints (\ref{eq:carton:c2b}, \ref{eq:carton:b2c}) define a rigid coupling
where selected cartons imply selected boxes and vice versa. This allows to
relax the box variables to be continuous while an optimal solution remains
binary.
In fact, the constraints imply the analytic relation
\begin{align}\label{eq:carton:box:candidates}
    y_b(z) &= 1-\prod_{(z_k,y_b)\in REL} (1-z_k)
\end{align}
which makes an optimal $y$ binary whenever $z$ is.
Moreover, the previous sub program \eqref{eq:sub:dual} can be used to
express the new one as $g(z)=f(y(z))$.
Consequently, the lower bound approximation \eqref{eq:sub:lowerbound}
requires the chain-rule for partial derivatives
\begin{equation} 
\begin{aligned}\label{eq:carton:dual:transform}
    g(z) & \geq s^i+ (w^i)'z \qquad s^i = f(y(z^i)) - (w^i)'z^i\\
    & w^i := \frac{dg}{dz}(z^i) = \left[\frac{df}{dy}(y(z^i))\right]\left[ \frac{dy}{dz}(z^i)\right] = \mu^i J^i \\
    & J^i_{bk} := \frac{dy_b}{dz_k}(z^i) = \prod_{\underset{k\neq l}{(z_l,y_b)\in REL}} (1-z^i_l)
\end{aligned}
\end{equation}
and builds on the optimal solution $f(y(z^i))$ of \eqref{eq:sub:dual}
corresponding to $y(z^i)$, with its optimal dual variables $\mu^i$,
and the Jaccobian $J^i_{bk}$ of \eqref{eq:carton:box:candidates} at $z^i$.
In other words, solving for $y$ is done ``in-between'' the master and sub problem
by a pre-/post-transformation of variables (c.f. \fig{fig:bender:summary}).

\usetikzlibrary{arrows,positioning,matrix,fit}
\tikzset{
    >=stealth',
    punkt/.style={
           rectangle,
           rounded corners,
           draw=black, very thick,
           text width=6.5em,
           minimum height=2em,
           text centered},
}

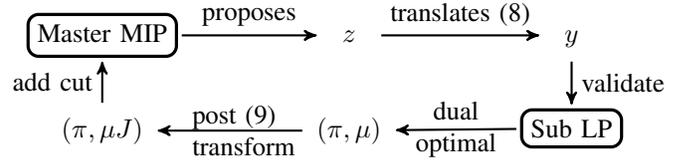
\begin{figure}[tb]
    \centering
    \begin{tikzpicture}[
            node distance=1cm, 
            box/.style={rounded corners, draw=black, rectangle, very thick},
            arr/.style={ ->, thick, shorten <=2pt, shorten >=2pt}
        ]
 \matrix (m) [matrix of math nodes,row sep=2em,column sep=5em,minimum width=2em] {
     \textrm{Master MIP}\vphantom{yM} & z\vphantom{yM} & y\vphantom{yM} \\
        (\pi, \mu J) & (\pi, \mu) & \textrm{Sub LP}\\
 };
        \node[box,fit=(m-1-1), inner sep=0pt]{};
        \node[box,fit=(m-2-3), inner sep=0pt]{};
        \draw (m-1-1) edge[arr] node[pos=0.5,above] {proposes} (m-1-2);
        \draw (m-1-2) edge[arr] node[pos=0.5,above] {translates \eqref{eq:carton:box:candidates}}(m-1-3);
        \draw (m-1-3) edge[arr] node[pos=0.5,right] {validate} (m-2-3);

        \draw (m-2-3) edge[arr] node[pos=0.5,text width=2cm,align=center] {dual\\optimal} (m-2-2);
        \draw (m-2-2) edge[arr] node[pos=0.5,text width=1cm,align=center] {post~\eqref{eq:carton:dual:transform}\\transform} (m-2-1);
        \draw (m-2-1) edge[arr] node[pos=0.5,left] {add cut} (m-1-1);
    \end{tikzpicture}%
    \caption{Benders iteration with intermediate variable pre-/post-transformations:
    (1) Solving the master program proposes carton candidates $z^i$.
    (2) They are translated by \eqref{eq:carton:box:candidates} to box candidates and validated with the sub program.
    (3) The dual optimal $\mu$ (with respect to boxes) are transformed back by 
        \eqref{eq:carton:dual:transform} to add new optimality cut (with respect to cartons) 
        in the master program.
    (4) This loop is repeated until convergence.
    }
    \label{fig:bender:summary}
\end{figure}

Last but not least, the remaining constraints \eqref{eq:box:fixed} have to be 
reformulated in terms of carton variables. This finally leads to the new master 
problem
\begin{equation}\label{eq:carton:master}
\begin{aligned}
    \min_{z,\theta}  \theta & \quad \text{subject to}\quad\\
    \theta &\geq \sum_p s^i + (w^i)'z \quad\forall i \\
    \sum_i z_k & = M \textrm{ and }
    \sum_{(z_k,y_b) \in REL} z_k &\geq 1 \quad \forall b \in \fixedBoxes \\
\end{aligned}
\end{equation}
where $s^i, w^i$ are computed via \eqref{eq:carton:dual:transform}.
In summary, 
the analytic transformation eliminates complexity in the master program with simple transformation of the dual optimal solution.
The new master problem depends only on $z, \theta$ and a
tiny number of constraints -- both favourable properties (c.f. \tab{tab:bender:summary}).

\begin{table}[tbp]
    \caption{Benders Decomposition Program Size Comparison}
\begin{center}
\begin{tabular}{|l|r|r|r|}
\hline
    Method & MIP Variables & MIP Constraints & LP Size \\
\hline
    No-Bender \ref{sec:problemstatement} & $\order{BP+K}$  & $\order{BP+K}$ & - \\
    Bender-x  \ref{sec:bender:packing}   & $\order{B+K+1}$ & $\order{B+K+I}$ & $\order{BP}$ \\
    Bender-xy \ref{sec:bender:carton}    & $\order{K+1}$   & $\order{1+I}$ & $\order{BP}$ \\
\hline
\end{tabular}\\[1mm]
    ($B\approx 70k$, $K\approx 40k$, $P\approx 110k$, Bender Iterations $I\leq 100$)
\end{center}
\label{tab:bender:summary}
\end{table}

\section{Efficient Benders Sub Program Implementation}\label{sec:subprogram}

Optimising the Benders sub program is essential for good performance. The linear
program is simpler than a mixed integer program, but its scale makes it
demanding, especially in terms of memory. The next section shows that the
Benders sub program has a closed-form, analytical solution that provides a
linear time algorithm. The following sections then exploits the structure of the
objective function to significantly reduce the memory requirements and
accelerate the evaluation.

\subsection{Analytical Benders sub problem}

The Benders sub problem \eqref{eq:sub:dual} is a linear program. This is much
simpler than a mixed integer program, but its challenging due to its scale ($P\times B$).
Note that the sub problem decomposes further since optimisation over packing units
are independent. This results in $P$ optimisations with each only $B$ variables.
In fact, those optimisations can be solved analytically
(assuming $C_{pb}=\infty$ when $F_{pb}=0$)
\begin{subequations}\label{eq:sub:analytical}
\begin{align}
    \pi_p    &= \min \{c_{pb} | F_{pb}=1, y_b=1\}\label{eq:sub:analytical:pi}\\
    \mu_{pb} &= \min (0, c_{pb} - F_{pb}\pi_p)\label{eq:sub:analytical:mu}
\end{align}
\end{subequations}
which the remainder of this section will show. As a result, solving the Benders
sub problem does not require an LP solver and runs in linear time.

Analysing one single packing unit $p$ in isolation yields
\begin{equation}
\begin{aligned}
    f^p_{primal}(y)  = \min_x \sum_{b} c_{pb} x_{b}
    & \quad\text{subject to}\\
    \forall b \quad x_{b} \leq y_b \textrm{ and }\quad & \sum_{b} F_{pb} x_{b} = 1
\end{aligned}
\end{equation}
The primal optimal solution is to choose \emph{the} best fitting box ($F_{pb}=1$) 
among the available ones ($y_b=1$). Note that this is what happens
on a daily basis when shipping a packing unit for real. This is a simple 
search for the smallest value in the subset as defined in \eqref{eq:sub:analytical:pi}. 
Solving the dual
\begin{equation}
\begin{aligned}
    f^p_{dual}(y) = \max_{\pi, \mu\leq 0} \pi + \sum_{b} \overbrace{y_b \mu_{b}}^{\leq 0}\\
\textrm{subject to}\quad \mu_{b} \leq c_{pb} - F_{pb}\pi \quad\forall b
\end{aligned}
\end{equation}
is more challenging. The key observation is that the second term is zero for
unselected boxes $y_b=0$ and never positive $\mu_{b} \leq \min(0,
c_{pb}-F_{pb}\pi_p) $ for selected boxes $y_b=1$. By choosing the best fitting
and available box \eqref{eq:sub:analytical:pi}, the second term is always zero
and $\pi$ matches the primal optimum. The $\mu_b$ corresponding $y_b=0$ are
not fully constrained, but choosing the largest possible value
\eqref{eq:sub:analytical:mu} yields the tightest lower bound in
\eqref{eq:sub:lowerbound}.

\subsection{Memory Efficient Implementation}\label{sec:subprogram:memory}

A naive implementation of \eqref{eq:sub:analytical} is to directly use the cost
$C_{pb}$. The problem is that explicitly ``materialising'' this cost matrix
requires a lot of memory. Without the analytical solution a solver would likely
use double precision floating point (as typical for numerical stability of a solver),
which would require about $60GB$ (for our dataset with $70k$ boxes and
$110k$ packing units). This is just the input data, let alone a solver's working data.

The analytical solution \eqref{eq:sub:analytical} allows for a more
memory efficient representation. We exploit the definition
$C_{pb}=\textrm{Volume}(b) - \textrm{Volume}(p)$ to compute the cost matrix 
on-the-fly based on the ``fitting matrix'' $F_{pb}$ and volumes. The volume
vectors are negligible $8\cdot(70k+110k)\approx 1.5 MB$, while the boolean matrix
$F$ can be bit-packed to require only $70k\cdot 110k\cdot 1\textrm{bit} \approx 1GB$.
This memory saving is paramount to work on typical commodity hardware with 
$64 GB$ of RAM.

\subsection{Runtime Efficient Implementation}\label{sec:subprogram:runtime}

Bit-packing also has positive impact on the runtime. First of all, loading
less data reduces I/O and is likely beneficial for CPU caches. Assuming a
64~bit architecture, bit-packing allows to load 64 values of $F$ with a single
CPU load instruction.
Dedicated CPU instructions of today's hardware also allow to efficiently
identify the non-zero bits $F_{pb}=1$. These are the only elements resulting in
a non-zero contribution to the overall result. 
The CPU instruction \emph{count trailing zeros} (\texttt{ctz})%
\footnote{\url{https://en.wikipedia.org/wiki/Find_first_set}} returns the
number of trailing 0-bits, starting at the least significant bit position. This
corresponds to an $F_{pb}=1$ which allows for computing $C_{pb}$ given the
indices. This bit can be cleared with simple bit operation after computing the
score contribution. The next 64~bits are loaded whenever the in-memory word
becomes zero. Note that this loop loads the full data matrix, but the actual
computation is only on non-zero elements (c.f. Algorithm~\ref{alg:dual}).

Furthermore, the boxes can be ordered such that the cost $C_{pb}$ is always
monotonic (in $b$). Consequently, finding $\pi_p$ does not always require a full scan
over all boxes, and all non-zero score contribution $\mu_{pb}$ for a given
packing unit $p$ have a smaller index than the best matching box. Hence, the
algorithm does not even process all boxes, and avoids the $\min(0, \cdot)$
in \eqref{eq:sub:analytical:mu}.

Last but not least, all packing units can be processed in parallel. We leverage
OpenMP \cite{openmp} in our implementation to get parallelisation with only
minor code annotations.

\newcommand*\BitAnd{\mathbin{\&}}
\newcommand*\BitOr{\mathbin{|}}
\newcommand*\ShiftLeft{\ll}
\newcommand*\ShiftRight{\gg}
\newcommand*\BitNeg{\ensuremath{\mathord{\sim}}}

\newcommand*\wordIndex{\textrm{wordIdx} }
\newcommand*\bitIndex{\textrm{bitIdx} }
\newcommand*\wordContent{\textrm{word} }

\begin{algorithm}[tb]
  \caption{FastDual($y$)}
\label{alg:dual}
\begin{algorithmic}
  \REQUIRE box indices $b$ ordered by volume $V(b)$
  \REQUIRE packing $F_{pb}$ packed into 64bit words
  \REQUIRE $\pi_p, \mu_{pb}$ initialized to zero
  \FORALL p
    \STATE bestBox := $\arg \min_b \{V_b | y_b=1, F_{pb}=1\}$
    \STATE $\pi_p$ += $V(bestBox) - V(p)$
    \STATE $\wordIndex, \bitIndex$ := divmod(bestBox, 64)
    \STATE \wordContent = $load(F, p, \wordIndex) \BitAnd (1 \ShiftLeft (63-\bitIndex) - 1)$
    \WHILE { $\wordIndex > 0\, \&\&\, \wordContent >0$}
      \WHILE { $\wordContent >0$ }
        \STATE $\bitIndex$ := $ctzll(\wordContent)$
        \STATE b := $64\cdot \wordIndex + (63-\bitIndex)$
        \STATE $\mu_{pb}$ += $V(b) - V(bestBox)$
        \STATE \wordContent := $\wordContent \BitAnd (\BitNeg (1 \ShiftLeft \bitIndex ))$
      \ENDWHILE
      \STATE \wordIndex := \wordIndex - 1
      \STATE \wordContent := $load(F, p, \wordIndex)$ (unless $\wordIndex <0$)
    \ENDWHILE
  \ENDFOR
\end{algorithmic}
\end{algorithm}

\section{KD-Tree for 3D Bin Packing Evaluation}\label{sec:kdtree}

The mixed integer program relies on the fitting matrix $F_{pb}$ which
indicates whether a packing unit $p$ fits into a box $b$. So far, we
assumed this matrix is given, while this section details its computation.

Computing a single entry $F_{pb}$ is known as the 3D bin packing problem
\cite{3dbinpacking}. The aim of this problem is to decide if a set of
rectangular items fits into a rectangular box. This problem is
NP-hard \cite{3dbinpacking-paper}. In general, it is also relevant how to fit
the items into the box, but we only require the boolean decision. We solve
a single $F_{pb}$ with an open source C\# implementation \cite{3dbinpacking-csharp}.

Our challenge is that the 3D bin packing problem needs to be solved not only
once, but for every pair of box and packing unit. In our experiment, that is
$B\times P\approx 7700M$ pairs. Even a single bin packing may be slow due to
the combinatorial growth of possibilities, especially when a packing unit
consists of many items. Hence, exhaustive evaluation is costly.

Many evaluations of the 3D bin packing can be deduced from others when making
the assumption that ``if a packing unit fits into a small box, it also fits into
larger ones''.%
\footnote{larger in all dimensions, not by volume}
We exploit this assumption using an adaptive KD-Tree partitioning.

\begin{figure}[tb]
    \centering
    \begin{tikzpicture}[
        scale=0.6,
        node distance=1cm, 
        box/.style={rounded corners, draw=black, rectangle, thick},
        anno/.style={rounded corners, fill=white},
        good/.style={rounded corners, fill=green,opacity=0.5},
        bad/.style={rounded corners, fill=red,opacity=0.5},
        recur/.style={rounded corners, fill=white,opacity=0.5}
    ]

    \coordinate (smallest) at (0,0);
    \coordinate (largest) at (10,10);
    \node[fit=(smallest)(largest),box]{};

    \draw (smallest) -- (largest);
    \coordinate (lower) at (6,6);
    \coordinate (upper) at (7,7);
    \fill[blue] (6.5, 6.5) circle (0.15) node[left]{\footnotesize binary search splitting};
        \node[fit=(smallest)(lower),box,bad] (TooSmall) {};
        \node[fit=(upper)(largest),box,good] (Fits) {};

    \coordinate (nw) at (0,10);
    \coordinate (se) at (6,7);
        \node[fit=(nw)(se),box,recur] (Splitting) {};

    \coordinate (ne) at (10,6);
    \coordinate (sw) at (7,0);
    \node[fit=(ne)(sw),box,recur] {};
    \coordinate (l1) at (7,0);
    \coordinate (u1) at (10,6);
    \draw (l1) -- node[pos=0.8](x){} (u1);
    \draw[fill=blue] (x) circle(0.15);
    \coordinate (l1a) at (9,4);
    \coordinate (u1b) at (10,5);
    \node[fit=(l1)(l1a),box,bad] {};
    \node[fit=(u1b)(u1),box,good] {};
    \coordinate (sw) at (7,5);
    \coordinate (ne) at (9,6);
    \node[fit=(sw)(ne),box,recur] {};
    \coordinate (sw) at (10,0);
    \coordinate (ne) at (10,4);
    \node[fit=(sw)(ne),box,recur] {};

    \coordinate (l1) at (7,5);
    \coordinate (u1) at (9,6);
    \draw (l1) -- node[pos=0.68](x){} (u1);
    \draw[fill=blue] (x) circle(0.15);
    \coordinate (l1a) at (8,5);
    \coordinate (u1b) at (9,6);
    \node[fit=(l1)(l1a),box,bad] {};
    \node[fit=(u1b)(u1),box,good] {};
    \coordinate (sw) at (7,6);
    \coordinate (ne) at (8,6);
    \node[fit=(sw)(ne),box,recur] {};
    \coordinate (sw) at (9,5);
    \coordinate (ne) at (9,5);
    \node[fit=(sw)(ne),box,recur] {};

    \foreach \x in {0,1,...,10} {
        \foreach \y in {0,1,...,10} {
            \fill[color=black] (\x,\y) circle (0.05);
        }
    }
    \node[anno,below right=7pt,text=red] at (TooSmall.north west) {\footnotesize boxes too small};
    \node[anno,below right=7pt,text=green!40!black] at (Fits.north west) {\footnotesize boxes fit};
    \node[anno,text width=25mm] at (Splitting) {\footnotesize boxes \emph{may} fit\\recursive partitioning};

    \end{tikzpicture}
    \caption{%
        Illustration of the KD-tree partitioning scheme -- in 2D instead of
				3D for simplicity.
        Every small black dot represents a box: the smallest box is bottom-left
        and the largest one is top-right.
        Bounding boxes (rectangle with rounded corner) represent a set of boxes.
        Binary search along the diagonal identifies a splitting point (blue circle)
        to recursively partition a set of boxes into quadrants (octants in 3D).
        All boxes in the upper-right (green) are large enough to fit a packing unit
        while boxes in the lower-left (red) are too small. All other (white)
        have to be partitioned recursively.
    }
\label{fig:kdtree}
\end{figure}
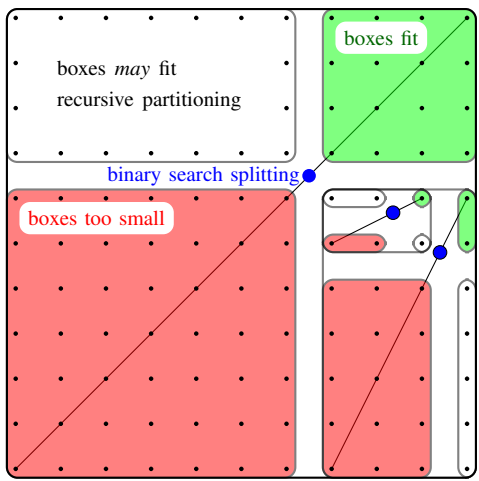

The KD-Tree based evaluation represents a set of boxes by their bounding box
(c.f. \fig{fig:kdtree}).
Every point within this bounding box represents a physical box. The evaluation
searches along the bounding box diagonal for the smallest box that still fits
all items. This point/box can be found with binary search. This partitions the
entire bounding box into octants, one of which contains strictly larger boxes,
one contains strictly smaller boxes, and six octants that have some dimensions
larger and others smaller. The binary decision $F_{pb}$ for the former two 
octants is clear, while the latter six octants need to be partitioned recursively.
When an octant is small enough (e.g. less 30 boxes), all contained boxes are directly
evaluated with the 3D bin packer.

\begin{figure}[tb]
    \centering
    \includegraphics[width=0.5\textwidth]{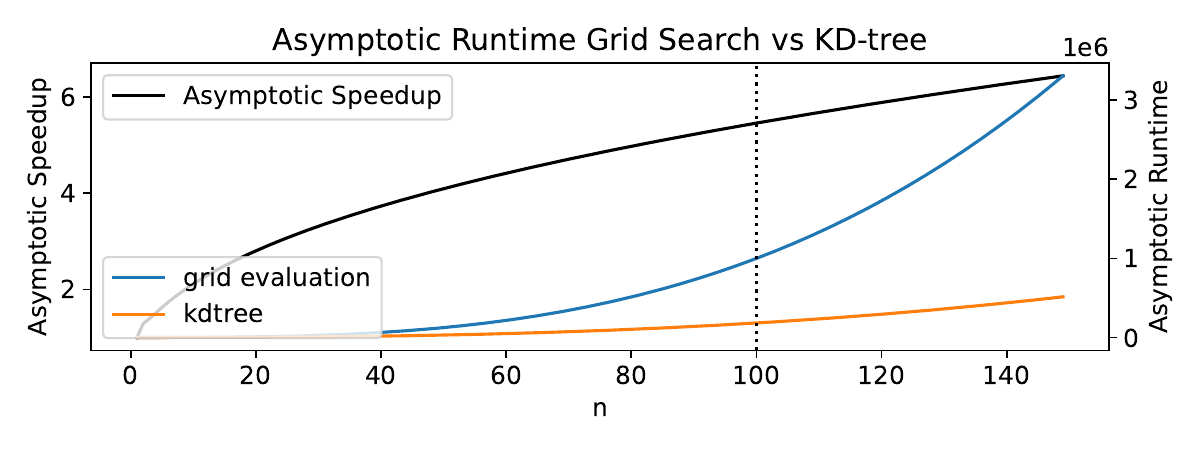}%
    \caption{%
        Asymptotic worst-case runtime of the KD-Tree evaluation $\order{n^{2.58}}$ 
        compared to exhaustive grid evaluation $\order{n^3}$. Discretisation of
        boxes of $1m$ to $1\cm$ precision is about $n\approx 100$ which suggests
        about a factor $5$ speedup.
    }%
\label{fig:kdtree:asymptotic}
\end{figure}

This KD-tree partitioning asymptotically requires only $\order{n^{2.58}}$ bin
packing evaluation when the bounding box is a cube with side length $n$. This
seems a small improvement over the exhaustive grid evaluation with $\order{n^3}$,
but results in about an order of magnitude speedup (c.f.
Figure~\ref{fig:kdtree:asymptotic}).

The runtime follows from telescoping
\footnote{
    See \url{https://www.wolframalpha.com/input?i=T(n)=log_2(n) \%2B 6T(n/2) and T(1)=1}
}
\begin{subequations}
\begin{align}
    T(n) = \overbrace{\log(n)}^{\textrm{binary search}} + 6 \overbrace{T(n/2)}^{\textrm{recursion}}\\
     = \log(n) + 6\log(n/2) + 6^2T(n/4) + \dots\\
     = \frac{31}{25} n^{\log_2(6)} -\frac{1}{5}\log_2(n) - \frac{6}{25} \\
     \approx \order{n^{\log 6}} \approx \order{n^{2.58496}}
\end{align}
\end{subequations}
with $\log_2$ the binary logarithm.
This is a worst case analysis when the splitting point is exactly half way,
i.e. all octants have equal volume and 6 of 8 have to be partitioned
recursively. If the split is off-center, the octants to be recurred on make up
less than $\frac{6}{8}$ of the bounding box volume, resulting in more
speedup.

\section{Results}\label{sec:results}

\subsection{Data Set and Evaluation Metrics}

All experiments were done with a real world data set. The data set consists
of $112'696$ packing units with $410'084$ items overall, i.e. about $4$
items per packing unit on average. The goal of the concrete optimisation
is to find $M=8$ cartons that minimise the empty space per packing. 
Additionally, three boxes (the largest possible and two with reduced
postage fee) where explicitly constrained to be part of the solution.

The largest box allowed by the Swiss postal service is $100 \times 60 \times 60 \cm$%
\footnote{\url{https://www.post.ch/de/pakete-versenden/pakete-schweiz/postpac-economy\#wichtig-zu-wissen}}.
Accounting for the carton thickness and the package label as a minimal size,
we consider all boxes between $15.5 \times 15.5 \times 10.5 \cm$ 
and $99.5\times59.5\times59.5\cm$
at steps of 
$1 \cm$. Furthermore, we assume $\textrm{length}\geq\textrm{width}\geq\textrm{height}$
to avoid symmetries, which results in a total of $B=71'790$ possible boxes.

The reported performance numbers are the objective value (empty space
in all packing units) normalised by the total volume of packing units,
which provides a more comparable score independent of the sample size.
The non-linear transformation $score/(1+score)$ then yields average empty
volume per box volume, which was the business' key performance indicator.

The following experiments were all done on a MacBook Pro (14-inch, 2021)
with ``M1 Max'' with 64 GB of RAM and 10 CPU cores. Our implementation uses
the open source \emph{Coin-or Branch and Cut (CBC)} solver \cite{coinor-cbc} 
and \emph{python-mip} \cite{python-mip} wrapper to solve mixed integer
programs.

\subsection{Evaluation of variable-height packaging optimisation}

\begin{figure}[tb]
    \centering
    \includegraphics[width=\columnwidth]{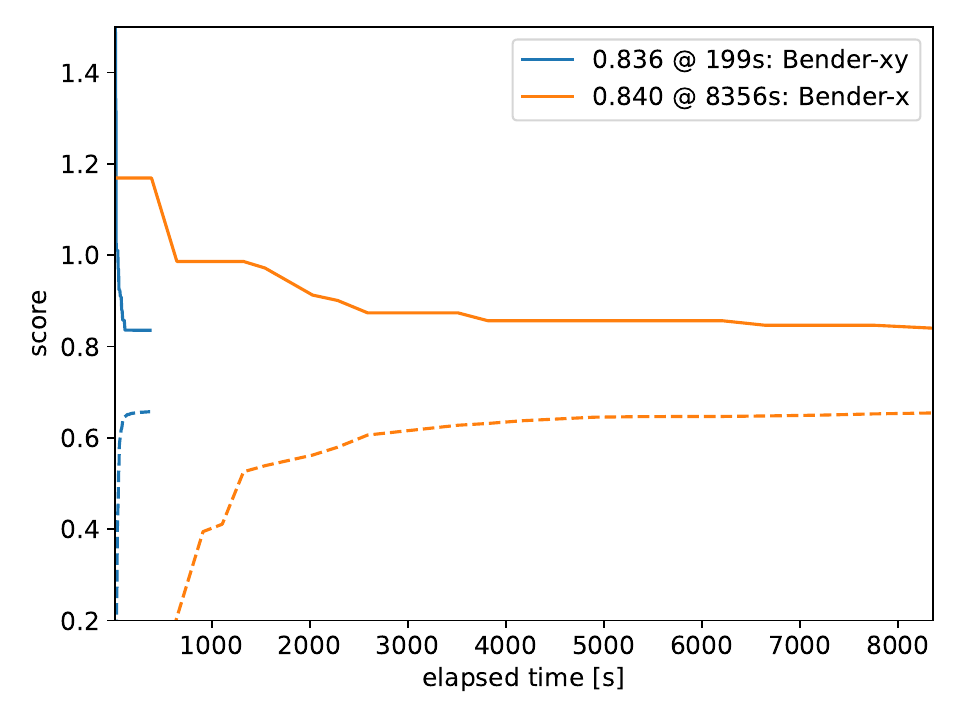}
    \caption{%
        Objective value and the theoretical lower bound for both
        variants of the Benders decomposition discussed in \secref{sec:bender}.
        Both variants yield comparable results and lower bounds. However,
        the variant \secref{sec:bender:carton} with fewer variables and constraints
        (blue) converges in less than $200s$ and is much faster, roughly by a factor 42.
    }
\label{fig:bender:result}
\end{figure}

\fig{fig:bender:result} reports the optimisation result along with the lower
bound provided by the Benders decomposition. The lower bound provides a
theoretical guarantee that no solution can be better and a solution is globally
optimal when the bound matches the score of the best found solution.
Unfortunately, the optimisation did not find the
global optimum and a relative gap of about 21\% remained.
This gap indicates that the optimal ``empty space'' is in the range of $[39.5\%\dots
45.5\%]$. The lower bound is a theoretical value, though, and the optimum is
more likely in between the two, making the gap smaller than it may appear.

We suspect the structure of the carton variable $z$ to be the reason for
stagnation. Recall that it forms a $43'256$ dimensional hypercube of which
$M=8$ cartons can be activated. Also, most cartons have 8 neighbours
(smaller/bigger in all 3 dimensions) resulting $8^8\approx 16M$ neighbouring
solutions. 
This likely requires many Benders optimality cuts to eliminate all near optimal
solutions, resulting in the stagnation.

\fig{fig:bender:result} also reports the benefit of moving most optimisation
variables into the Benders sub program. The chart shows that having only the
carton variables $z$ in the master program outperforms the alternative
significantly. Keeping the box variables $y$ in the master program increases
the number of optimisation variables and adds many complicating constraints, which
massively slows down the convergence measured in terms of time.

In summary, simplifying the master program as much as possible is key.
Moving the box variables $y$ to the sub program is cheap, it only requires 
a simple transformation from and to the carton variable, i.e.
\eqref{eq:carton:box:candidates} and \eqref{eq:carton:dual:transform}.
This results in computational savings from $8'356s$ to $199s$, i.e. a speedup
of two orders of magnitude, without loss of accuracy.

\subsection{KD-Tree versus exhaustive grid bin packing evaluation}

The mixed integer optimisation requires the bin packing matrix $F_{pb}$ as
input. This boolean matrix is computed prior to the mixed integer program and
is currently the key performance bottleneck. \tab{fig:kdtree:runtime} compares
the exhaustive grid evaluation to our adaptive KD-Tree evaluation. 

\begin{table}[htbp]
    \caption{%
        Comparing exhaustive grid and KD-Tree evaluation for $F_{pb}$
    }
\begin{center}
\begin{tabular}{|c|c|c|}
\hline
    &Exhaustive Grid & Adaptive KDTree\\
\hline
    CPU-time  & 535.76h & 18.40h \\
    Elapsed (10 core) & 74h17m & 2h59m\\
\hline
    Empty Volume & 45.3\% & 45.5\% \\
    Best Score & 0.829 & 0.836 \\
    Lower Bound & 0.655 & 0.657 \\
\hline
\end{tabular}
\end{center}
\label{fig:kdtree:runtime}
        The KD-Tree provides a speedup of 29 in terms of CPU-time with
        comparable result.
\end{table}

The runtime was measured with UNIX \texttt{time} utility showing the CPU-time (as if
run on a single core) and overall elapsed time. The latter is not a good
performance measure, as the evaluations were run as a background process with
low ``nice'' value. But it provides a rough estimate of how long it runs, while
CPU-time is the more robust performance metric.

Overall, the KD-Tree runs only about 3 hours instead of multiple days. Looking
at the CPU-time, we observe more than an order of magnitude speedup. The
evaluation of packing matrix $F$ is the most expensive part and having a
speedup of 29 is very important for end-to-end runtime.

\subsection{Optimised Benders sub program performance}

\begin{table}[tbp]
    \caption{%
        Performance comparison of efficient dual implementation
    }
\begin{center}
\begin{tabular}{|r|r||r|r|r|}
\hline
    P & Cores & Fast-C & Naive-C  & Python \\
\hline
 10000 &  1 & 0.0302s  &1.2065s  &7.1931s\\ 
 10000 & 10 & 0.0069s  &0.1769s  & -  \\
\hline
112696 & 10 & 0.0618s  &1.9398s  &  out-of-mem \\
\hline
\end{tabular}
\end{center}
\label{fig:dual:fast}
        (upper) Comparision with a non-bit optimised Python implementation
        is done on a $P=10k$ subset to avoid unnecessary memory pressure,
        and timing is averaged over 10 runs.
        A naive bit-packed C implementation (without ctz) yields a 6 fold speed over Python, while
        exploiting bit-optimised CPU instructions yields another 40 fold speedup.
        OpenMP parallelisation yields another 4 fold improvement on a 10 core M1 Max.
        In total, the bit-optimised, multi-core implementation is 3 orders of magnitude
        faster $7.2/0.007=1029$.
        (lower) Performance on full data set (averaged over 100 runs) runs
        in less than $70ms$, making the sub program essentially ``cost-free''.
\end{table}

\tab{fig:dual:fast} reports the performance of the Bender sub program optimisation
(\secref{sec:subprogram}).
Overall, the implementation yields an improvement of 3 orders of magnitude. The speedup
builds on three complementary aspects:
\begin{itemize}
    \item 6 fold speedup: bit-packed data representation using C instead of Python implementation
    \item 40 fold speedup: loading 64 bit words and processing only non-zero bits
        using dedicated $ctzll$ CPU instruction
    \item 4 fold speedup: parallelisation using OpenMP\cite{openmp}
\end{itemize}

The bit-packed representation is absolutely crucial, not only for performance,
but to reduce memory usage by a factor 64, i.e. about $1GB$ instead of $64GB$.
This is key to run the optimisation on commodity hardware.
Moreover, using dedicated CPU instructions shows significantly more speedup over
``just'' parallelising on multiple cores. In fact, processing only non-zero bits
exploits the sparsity of $F_{pb}$ which explains the tremendous speedup.%
\footnote{ A sparse matrix representation would probably be less effective due
to more irregular memory access and computation overhead. }

All optimisations combined make solving the Benders sub program almost cost-free.
This is a major performance improvement over a na\"ive implementation and make
experiments on an ultra portable laptop possible.

Also recall that the memory usage and runtime of the sub program are linear in
the number of packing units $P$ and that the master program's runtime does not depend on $P$.
This is beneficial as using a larger data sample will have little impact
on the overall optimisation and allows for experimenting with larger data sets in the
future.

\section{Future Work \& Conclusion}\label{sec:conclusion}

\subsection{Future Work}

The presented carton optimisation model focuses on reducing the overall empty
space per packing unit. There are however additional metrics worth looking at.

While crease lines allow for making multiple (smaller) boxes out of a carton,
it still requires the packer to make small cuts into the carton. This costs
additional time compared to only folding the carton into the largest box. This
cost could be incorporating into the objective in the future.

Overall, the performance optimisation presented in this work provides a solid
basis for further investigation of the above topics. Concerning the Benders
decomposition it might be interesting to build on solver callbacks for further
improvements \cite{benders-then-and-now}. This allows to interleave solving the
sub program and the master program to avoid recomputing the master program over
and over again.

\subsection{Conclusion}

This paper tackles the problem of finding optimal variable-height transport
packaging. This is a challenge faced by online retailers when optimising their
shipment process.

We present a mixed integer formulation that models cartons with crease lines.
It aims for optimal boxes to ship goods with the least amount of empty space.
Minimising empty space is crucial to save on filler, reduce packing time,
and improve customer satisfactions.

The key contributions of this paper are problem-specific optimisations to solve
that mixed integer program efficiently. The resulting model had over seven
billion variables initially and solving it requires optimisations on many
different levels.

\emph{Mathematically}, the Benders decomposition is what makes all possible. Our
contribution is how to incorporate the crease lines into the Bender
decomposition to achieve a 42-fold speedup, with simple
transformation of variables.

\emph{Algorithmically}, our KD-Tree based evaluation of the 3D bin packing showed a
29-fold speedup. Bin-packing is currently the bottleneck and this optimisation 
reduces the overall optimisation to only about three hours on a laptop.

\emph{Computationally}, we show an astonishing speedup of 3 orders of magnitude for
evaluating the Benders sub program. This optimisation built on the analytical
solution of the linear program and used algorithmic improvement as well as
specialised CPU instruction of modern hardware. This makes the sub program
almost free of charge and allows for larger data sample in the future.

\bibliographystyle{IEEEtrans} 
\bibliography{IEEEabrv,references} 

\begin{thebibliography}{10}
\providecommand{\url}[1]{#1}
\csname url@samestyle\endcsname
\providecommand{\newblock}{\relax}
\providecommand{\bibinfo}[2]{#2}
\providecommand{\BIBentrySTDinterwordspacing}{\spaceskip=0pt\relax}
\providecommand{\BIBentryALTinterwordstretchfactor}{4}
\providecommand{\BIBentryALTinterwordspacing}{\spaceskip=\fontdimen2\font plus
\BIBentryALTinterwordstretchfactor\fontdimen3\font minus
  \fontdimen4\font\relax}
\providecommand{\BIBforeignlanguage}[2]{{%
\expandafter\ifx\csname l@#1\endcsname\relax
\typeout{** WARNING: IEEEtranS.bst: No hyphenation pattern has been}%
\typeout{** loaded for the language `#1'. Using the pattern for}%
\typeout{** the default language instead.}%
\else
\language=\csname l@#1\endcsname
\fi
#2}}
\providecommand{\BIBdecl}{\relax}
\BIBdecl

\bibitem{3dbinpacking-paper}
\BIBentryALTinterwordspacing
E.~Baltacioglu, ``The distributer's three-dimensional pallet-packing problem: A
  human intelligence-based heuristic approach,'' Master's thesis, Air Force
  Institute of Technology, 2001. [Online]. Available:
  \url{https://scholar.afit.edu/etd/4563}
\BIBentrySTDinterwordspacing

\bibitem{benders}
J.~F. Benders, ``Partitioning procedures for solving mixed-variables
  programming problems ‘,'' \emph{Numerische mathematik}, vol.~4, no.~1, pp.
  238--252, 1962.

\bibitem{bendersCplex}
P.~Bonami, D.~Salvagnin, and A.~Tramontani, ``Implementing automatic benders
  decomposition in a modern mip solver,'' in \emph{Integer Programming and
  Combinatorial Optimization}, D.~Bienstock and G.~Zambelli, Eds.\hskip 1em
  plus 0.5em minus 0.4em\relax Cham: Springer International Publishing, 2020,
  pp. 78--90.

\bibitem{openmp}
R.~Chandra, L.~Dagum, D.~Kohr, R.~Menon, D.~Maydan, and J.~McDonald,
  \emph{Parallel programming in OpenMP}.\hskip 1em plus 0.5em minus 0.4em\relax
  Morgan kaufmann, 2001.

\bibitem{3dbinpacking-csharp}
D.~Chapman, ``3dcontainerpacking,''
  \url{https://github.com/davidmchapman/3DContainerPacking}, 2021.

\bibitem{Fischetti17}
\BIBentryALTinterwordspacing
M.~Fischetti, I.~Ljubić, and M.~Sinnl, ``{Redesigning Benders Decomposition
  for Large-Scale Facility Location},'' \emph{Management Science}, vol.~63,
  no.~7, pp. 2146--2162, July 2017. [Online]. Available:
  \url{https://ideas.repec.org/a/inm/ormnsc/v63y2017i7p2146-2162.html}
\BIBentrySTDinterwordspacing

\bibitem{coinor-cbc}
\BIBentryALTinterwordspacing
J.~Forrest, T.~Ralphs, H.~G. Santos, S.~Vigerske, J.~Forrest, L.~Hafer,
  B.~Kristjansson, jpfasano, EdwinStraver, M.~Lubin, rlougee, jpgoncal1,
  Jan-Willem, h-i gassmann, S.~Brito, Cristina, M.~Saltzman, tosttost,
  B.~Pitrus, F.~MATSUSHIMA, and to~st, ``coin-or/cbc: Release
  releases/2.10.8,'' May 2022. [Online]. Available:
  \url{https://doi.org/10.5281/zenodo.6522795}
\BIBentrySTDinterwordspacing

\bibitem{Gurumoorthy2020}
\BIBentryALTinterwordspacing
K.~S. Gurumoorthy, S.~Sanyal, and V.~Chaoji, ``Think out of the package:
  Recommending package types for e-commerce shipments,'' in \emph{ECML-PKDD
  2020}, 2020. [Online]. Available:
  \url{https://www.amazon.science/publications/think-out-of-the-package-recommending-package-types-for-e-commerce-shipments}
\BIBentrySTDinterwordspacing

\bibitem{3dbinpacking}
S.~Martello, D.~Pisinger, and D.~Vigo, ``The three-dimensional bin packing
  problem,'' \emph{Operations Research}, vol.~48, 02 1998.

\bibitem{ortec}
B.~Miller, ``Analyzing dimensional weight pricing and load optimization,''
  Ortec, Tech. Rep., 2015.

\bibitem{bendersSurvey}
R.~Rahmaniani, T.~G. Crainic, M.~Gendreau, and W.~Rei, ``The benders
  decomposition algorithm: A literature review,'' \emph{European Journal of
  Operational Research}, vol. 259, no.~3, pp. 801--817, 2017.

\bibitem{benders-then-and-now}
\BIBentryALTinterwordspacing
P.~A. Rubin, ``Benders decomposition then and now,'' Oct. 2011. [Online].
  Available:
  \url{https://orinanobworld.blogspot.com/2011/10/benders-decomposition-then-and-now.html}
\BIBentrySTDinterwordspacing

\bibitem{python-mip}
\BIBentryALTinterwordspacing
T.~A.~M. Toffolo and H.~G. Santos, ``python-mip: Release releases/2.10.8,''
  Dec. 2022. [Online]. Available: \url{https://python-mip.com}
\BIBentrySTDinterwordspacing

\bibitem{adaptivemulticut}
\BIBentryALTinterwordspacing
S.~Trukhanov, L.~Ntaimo, and A.~Schaefer, ``Adaptive multicut aggregation for
  two-stage stochastic linear programs with recourse,'' \emph{European Journal
  of Operational Research}, vol. 206, no.~2, pp. 395--406, 2010. [Online].
  Available:
  \url{https://www.sciencedirect.com/science/article/pii/S0377221710001566}
\BIBentrySTDinterwordspacing

\bibitem{VANLOON2015478}
\BIBentryALTinterwordspacing
P.~{van Loon}, L.~Deketele, J.~Dewaele, A.~McKinnon, and C.~Rutherford, ``A
  comparative analysis of carbon emissions from online retailing of fast moving
  consumer goods,'' \emph{Journal of Cleaner Production}, vol. 106, pp.
  478--486, 2015. [Online]. Available:
  \url{https://www.sciencedirect.com/science/article/pii/S0959652614006489}
\BIBentrySTDinterwordspacing

\bibitem{buycom}
C.~Weber, C.~Hendrickson, P.~Jaramillo, S.~Matthews, A.~Nagengast, and
  R.~Nealer, ``Life cycle comparison of traditional retail and e-commerce
  logistics for electronic products: A case study of buy.com,'' Carnegie Mellon
  University, Tech. Rep., 2008.

\end{thebibliography}

\end{document}